\documentclass{amsart}
\usepackage{latexsym}
\usepackage[T1]{fontenc}
\usepackage{amsfonts}
\usepackage{amscd}
\usepackage{amssymb}
\usepackage{mathrsfs}
\usepackage{amsmath}
\usepackage{amsthm}
\newtheorem{thm}{Theorem}[section]
\newtheorem{cor}[thm]{Corollary}

\newtheorem{pro}[thm]{Proposition}
\theoremstyle{definition}
\newtheorem{defn}[thm]{Definition}
\newtheorem{ex}[thm]{Example}
\newtheorem{rmk}[thm]{Remark}
\DeclareMathOperator{\Per}{Per}
\DeclareMathOperator{\supp}{supp}

\DeclareMathOperator{\sep}{Sep}

\DeclareMathOperator{\Aut}{Aut}
\DeclareMathOperator{\Homeo}{Homeo}
\DeclareMathOperator{\id}{id}

\begin{document}
\author{Christian Svensson}
\address{Mathematical Institute, Leiden University,
P.O. Box 9512, 2300 RA Leiden, The Netherlands, and Centre for
Mathematical Sciences, Lund University, Box 118, SE-221 00 Lund,
Sweden} \email{chriss@math.leidenuniv.nl}
\author{Sergei Silvestrov}
\address{Centre for Mathematical Sciences,
Lund University, Box 118, SE-221 00 Lund, Sweden}
\email{Sergei.Silvestrov@math.lth.se}
\author{Marcel de Jeu}
\address{Mathematical Institute,
Leiden University, P.O. Box 9512, 2300 RA Leiden, The Netherlands}
\email{mdejeu@math.leidenuniv.nl}

\title{DYNAMICAL SYSTEMS AND COMMUTANTS IN CROSSED PRODUCTS}
\noindent 
\noindent  \noindent

\maketitle

\begin{abstract}
In this paper we describe the commutant of an arbitrary subalgebra $A$ of the
algebra of functions on a set $X$ in a crossed product of $A$ with the
integers, where the latter act on $A$ by a composition automorphism defined via
a bijection of $X$. The resulting conditions which are necessary and sufficient
for $A$ to be maximal abelian in the crossed product are subsequently applied
to situations where these conditions can be shown to be equivalent to a
condition in topological dynamics. As a further step, using the Gelfand
transform we obtain for a commutative completely regular semi-simple Banach
algebra a topological dynamical condition on its character space which is
equivalent to the algebra being maximal abelian in a crossed product with the
integers.
\end{abstract}
\section{Introduction}
The description of commutative subalgebras in non-commutative algebras and
their properties is an important direction of investigation for any class of
non-commutative algebras and rings, because it allows one to relate
representation theory, non-commutative properties, graded structures, ideals
and subalgebras, homological and other properties of non-commutative algebras
to spectral theory, duality, algebraic geometry and topology naturally
associated with the commutative subalgebras. In representation theory, for
example, one of the keys for the construction and classification of
representations is the method of induced representations. The underlying
structures behind this method are the semi-direct products or crossed products
of rings and algebras by various actions. When a non-commutative algebra is
given, one looks for a subalgebra such that its representations can be studied
and classified more easily, and such that the whole algebra can be decomposed
as a crossed product of this subalgebra by a suitable action. Then the
representations for the subalgebra are extended to representations of the whole
algebra using the action and its properties. A description of representations
is most tractable for commutative subalgebras as being, via the spectral theory
and duality, directly connected to algebraic geometry, topology or measure
theory.

If one has found a way to present a non-commutative algebra as a crossed
product of a commutative subalgebra by some action on it of the elements from
outside the subalgebra, then it is important to know whether this subalgebra is
maximal abelian or, if not, to find a maximal abelian subalgebra containing the
given subalgebra, since if the selected subalgebra is not maximal abelian, then
the action will not be entirely responsible for the non-commutative part as one
would hope, but will also have the commutative trivial part taking care of the
elements commuting with everything in the chosen commutative subalgebra. This
maximality of a commutative subalgebra and related properties of the action are
intimately related to the description and classifications of representations of
the non-commutative algebra.

Little is known in general about connections between properties of the
commutative subalgebras of crossed product algebras and properties of dynamical
systems that are in many situations naturally associated with the construction.
A remarkable result in this direction is known, however, in the context of
crossed product $C^*$-algebras. When the algebra is described as the crossed
product $C^*$-algebra $C(X) \rtimes_{\alpha} \mathbb{Z}$ of the algebra of
continuous functions on a compact Hausdorff space $X$ by an action of
$\mathbb{Z}$ via the composition automorphism associated with a homeomorphism
$\sigma$ of $X$, it is known that $C(X)$ sits inside the $C^*$-crossed product
as a maximal abelian subalgebra if and only if for every positive integer $n$
the set of points in $X$ having period $n$ under iterations of $\sigma$ has no
interior points  \cite{LiBing-Ren}, \cite{Tom1}, \cite{Tom2}, \cite{Tom3}, \cite{Zeller-Meier}. By the category theorem, this condition is
equivalent to the action of $\mathbb Z$ on $X$ being topologically free in the
sense that the non-periodic points of $\sigma$ are dense in $X$. This result on
the interplay between the topological dynamics of the action on the one hand,
and the algebraic property of the commutative subalgebra in the crossed product
being maximal abelian on the other hand, provided the main motivation and
starting point for our work.

In this article, we bring such interplay results into a more general algebraic
and set theoretical context of a crossed product $A\rtimes_{\alpha}\mathbb{Z}$
of an arbitrary subalgebra $A$ of the algebra $\mathbb{C}^X$ of functions on a
set $X$ (under the usual pointwise operations) by $\mathbb{Z}$, where the
latter acts on $A$ by a composition automorphism defined via a bijection of
$X$. In this general algebraic set theoretical framework the essence of the
matter is revealed. Topological notions are not available here and thus the
condition of freeness of the dynamics as described above is not applicable, so
that it has to be generalized in a proper way in order to be equivalent to the
maximal commutativity of $A$. We provide such a generalization. In fact, we
describe explicitly the (unique) maximal abelian subalgebra containing $A$
(Theorem \ref{commdesc}), and then the general result, stating equivalence of
maximal commutativity of $A$ in the crossed product and the desired
generalization of topological freeness of the action, follows immediately
(Theorem \ref{sepdom}). It involves separation properties of $A$ with respect
to the space $X$ and the action.

The general set theoretical framework allows one to investigate the relation
between the maximality of the commutative subalgebra in the crossed product on
the one hand, and the properties of the action on the space on the other hand,
for arbitrary choices of the set $X$, the subalgebra $A$ and the action,
different from the previously cited classical choice of continuous functions
$C(X)$ on a compact Hausdorff topological space $X$. As a rather general
application of our results we obtain that, for a Baire topological space $X$
and a subalgebra $A$ of $C(X)$ satisfying a mild separation condition, the
property of $A$ being maximal abelian in the crossed product is equivalent to
the action being topologically free in the sense above, i.e., to the set of
non-periodic points being dense in $X$ (Theorem \ref{bairemax}). This applies
in particular when $X$ is a compact Hausdorff space, so that a result analogous
to that for the crossed product $C^*$-algebra $C(X) \rtimes_{\alpha}
\mathbb{Z}$ is obtained. We also demonstrate that, for a general topological
space $X$ and a subalgebra $A$ of $C(X)$ satisfying a less common separation
condition, the subalgebra $A$ is maximal abelian if and only if $\sigma$ is not
of finite order, which is a much less restrictive condition than topological
freeness (Theorem \ref{domunmax}). Examples of this situation are provided by
crossed products of subalgebras of holomorphic functions on a connected complex
manifolds by biholomorphic actions (Corollary \ref{comfd}). It is interesting
to note that these two results, Theorem \ref{bairemax} and Theorem
\ref{domunmax}, have no non-trivial situations in common (Remark~\ref{disrmk}).

In the motivating Section \ref{sec:motivation}, we illustrate by several
examples that, generally speaking, topological freeness and the property of the
subalgebra $A$ being maximal abelian in the crossed product are unrelated. In
Section \ref{sec:dualsystem}, we consider the crossed product of a commutative
completely regular semi-simple Banach algebra $A$ by the action of an
automorphism. Since the automorphism of $A$ induces a homeomorphism of $\Delta
(A)$, the set of all characters on $A$ endowed with the Gelfand topology, such
a crossed product is via the Gelfand transform canonically isomorphic (by
semi-simplicity) to the crossed product of the subalgebra $\widehat{A}$ of the
algebra $C_0 (\Delta(A))$ by the induced homeomorphism, where $C_0 (\Delta(A))$ denotes the algebra of continuous functions on $\Delta(A)$ that vanish at infinity. Since $A$ is also
assumed to be completely regular, the general result in Theorem \ref{bairemax}
can then be applied to this isomorphic crossed product. We thus prove that,
when $A$ is a commutative completely regular semi-simple Banach algebra, it is
maximal abelian in the crossed product if and only if the associated dynamical
system on the Gelfand dual $\Delta (A)$ is topologically free in the sense that
its non-periodic points are dense in the topological space $\Delta (A)$
(Theorem \ref{compregmax}). It is then possible to understand the algebraic
properties in the examples in Section \ref{sec:motivation} as dynamical
properties after all, by now looking at the ``correct'' dynamical system on
$\Delta(A)$ instead of the originally given dynamical system. In these
examples, this corresponds to removing parts of the space or enlarging it.

In Section \ref{sec:loccompgr}, we apply Theorem \ref{compregmax} in the case
of an automorphism of $L_1(G)$, the commutative completely regular semi-simple
Banach algebra consisting of equivalence classes of complex valued Borel
measurable functions on a locally compact abelian group $G$ that are integrable
with respect to a Haar measure, with multiplication given by convolution. We
use a result saying that every automorphism on $L_1(G)$ is induced by a piecewise
affine homeomorphism on the dual group of $G$ (Theorem \ref{autohomeo}) to
prove that, for $G$ with connected dual, $L_1(G)$ is maximal abelian in the
crossed product if and only if the given automorphism of $L_1(G)$ is not of
finite order (Theorem \ref{conndual}). We also provide an example showing that
this equivalence may fail when the dual of $G$ is not connected.

Finally, in Section \ref{sec:generatorscom}, we provide another application of
the commutant description in Theorem \ref{commdesc} and of the isomorphic
crossed product on the character space, showing that, for a semi-simple
commutative Banach algebra $A$, the commutant of $A$ in the crossed product is
finitely generated as an algebra over $\mathbb{C}$ if and only if $A$ has
finite dimension as a vector space over $\mathbb{C}$.
\section{Crossed products associated with automorphisms}
\subsection{Definition}
Let $A$ be an associative commutative $\mathbb{C}$-algebra and let $\Psi : A \rightarrow A$ be an automorphism.
Consider the set
\[A \rtimes_{\Psi} \mathbb{Z} = \{f: \mathbb{Z} \rightarrow A \,|\, f(n) = 0 \,\,\textup{except for a finite number of}\, \,n\}.\]
We endow it with the structure of an associative $\mathbb{C}$-algebra by defining scalar multiplication and addition as the usual pointwise operations.
Multiplication is defined by \emph{twisted convolution}, $*$, as follows;
\[(f*g) (n) = \sum_{k \in \mathbb{Z}} f(k) \cdot \Psi^k (g(n-k)),\]
where $\Psi^k$ denotes the $k$-fold composition of $\Psi$ with itself.
It is trivially verified that $A \rtimes_{\sigma} \mathbb{Z}$ \emph{is} an associative $\mathbb{C}$-algebra under these operations. We call it the \emph{crossed product of $A$ and $\mathbb{Z}$ under $\Psi$}.

A useful way of working with $A \rtimes_{\Psi} \mathbb{Z}$ is to write elements $f, g \in A \rtimes_{\Psi} \mathbb{Z}$ in the form
$f = \sum_{n \in \mathbb{Z}} f_n \delta^n, g = \sum_{m \in \mathbb{Z}} g_m \delta^n$, where $f_n = f(n), g_m = g(m)$, addition and scalar multiplication are canonically defined, and multiplication is determined by $(f_n \delta^n)*(g_m \delta^m) = f_n \cdot \Psi^n (g_m) \delta^{n+m}$, where $n,m \in \mathbb{Z}$ and $f_n, g_m \in A$ are arbitrary. Using this notation, we may think of the crossed product as a complex Laurent polynomial algebra in one variable (having $\delta$ as its indeterminate) over $A$ with twisted multiplication.
\subsection{A maximal abelian subalgebra of $A \rtimes_{\Psi} \mathbf{Z}$}
Clearly one may canonically view $A$ as an abelian subalgebra of $A \rtimes_{\Psi} \mathbb{Z}$, namely as $\{f_0 \delta^0 \,|\, f_0 \in A\}$. Here we prove that the commutant of $A$ in $A \rtimes_{\Psi} \mathbb{Z}$, denoted by $A'$, is commutative, and thus there exists a \emph{unique} maximal abelian subalgebra of $A \rtimes_{\Psi} \mathbb{Z}$ containing $A$, namely $A'$.

Note that we can write the product of two elements in  $A \rtimes_{\Psi} \mathbb{Z}$ as follows;
\[f * g = (\sum_{n \in \mathbb{Z}} f_n \delta^n) * (\sum_{m \in \mathbb{Z}} g_m \delta^m) = \sum_{r, n \in \mathbb{Z}} f_n \cdot \Psi^n (g_{r-n}) \delta^r.\]
We see that the precise condition for commutation of two such elements $f$ and $g$ is
\begin{equation}\forall r: \sum_{n \in \mathbb{Z}} f_n \cdot \Psi^n (g_{r-n}) = \sum_{m \in \mathbb{Z}} g_m \cdot \Psi^m (f_{r-m}). \label{ekv}\end{equation}
\begin{pro}\label{commabel}
The commutant $A'$ is abelian, and thus it is the unique maximal abelian subalgebra containing $A$.
\end{pro}
\begin{proof}
Let $f = \sum_{n \in \mathbb{Z}} f_n \delta^n, g = \sum_{m \in \mathbb{Z}} g_m \delta^m \in A \rtimes_{\Psi} \mathbb{Z}$. We shall verify that if $f, g \in A'$, then $f$ and $g$ commute.
Using (\ref{ekv}) we see that membership of $f$ and $g$ respectively in $A'$ is equivalent to the equalities
\begin{align}&\forall h \in A, \forall n \in \mathbb{Z} : f_n \cdot \Psi^{n} (h) = f_n \cdot h, \label{ekv2}\\
&\forall h \in A, \forall m \in \mathbb{Z} : g_m \cdot \Psi^m (h) = g_m \cdot h. \label{ekv3}\end{align}
Insertion of (\ref{ekv2}) and (\ref{ekv3}) into (\ref{ekv}) we realize that the precise condition for commutation of such $f$ and $g$ can be rewritten as
\[\forall r: \sum_{n \in \mathbb{Z}} f_n \cdot g_{r-n} = \sum_{m \in \mathbb{Z}} g_m \cdot f_{r-m}.\]
This clearly holds, and thus $f$ and $g$ commute. From this it follows immediately that $A'$ is the unique maximal abelian subalgebra containing $A$.
\end{proof}
\section{Automorphisms induced by bijections}\label{auto}
Fix a non-empty set $X$, a bijection $\sigma : X \rightarrow X$, and an algebra of functions $A \subseteq \mathbb{C}^X$ that is invariant under $\sigma$ and $\sigma^{-1}$, i.e., such that if $h \in A$, then $h \circ \sigma \in A$ and $h \circ \sigma^{-1} \in A$.
Then $(X, \sigma)$ is a discrete dynamical system (the action of $n \in \mathbb{Z}$ on $x \in X$ is given by $n: x \mapsto \sigma^n (x)$) and $\sigma$ induces an automorphism $\widetilde{\sigma} : A \rightarrow A$ defined by $\widetilde{\sigma} (f) = f \circ \sigma^{-1}$ by which $\mathbb{Z}$ acts on $A$ via iterations.

In this section we will consider the crossed product $A \rtimes_{\widetilde{\sigma}} \mathbb{Z}$ for the above setup, and explicitly describe the commutant, $A'$, of $A$ and then center, $Z(A \rtimes_{\widetilde{\sigma}} \mathbb{Z})$. Furthermore, we will investigate equivalences between properties of non-periodic points of the system $(X, \sigma)$, and properties of $A'$.
First we make a few definitions.
\begin{defn}\label{sepper}
For any nonzero $n \in \mathbb{Z}$ we set
\begin{align*}
\sep_A^n (X) &= \{x \in X | \exists h \in A :  h(x) \neq \widetilde{\sigma}^n (h) (x)\}, \\
\Per_A^n (X) &= \{x \in X | \forall h \in A : h(x) = \widetilde{\sigma}^n (h) (x)\},\\
\sep^n (X) &= \{x \in X |x \neq \sigma^{n} (x))\}, \\
\Per^n (X) &= \{x \in X |x = \sigma^{n} (x)\}.
 \end{align*}
Furthermore, let
\begin{align*}
\Per_A^{\infty} (X) &= \bigcap_{n \in \mathbb{Z} \setminus \{0\}} \sep_A^n (X),\\
\Per^{\infty}(X) &= \bigcap_{n \in \mathbb{Z} \setminus \{0\}} \sep^n (X).\\
\end{align*}
Finally, for $f \in A$, put
\begin{align*}\supp(f) &= \{x \in X \, | \, f(x) \neq 0\}.
\end{align*}
\end{defn}
It is easy to check that all these sets, except for $\supp(f)$, are $\mathbb{Z}$-invariant and that if $A$ separates the points of $X$, then $\sep_A^n (X) = \sep^n (X)$ and $\Per_A^n (X) = \Per^n(X)$.
Note also that $X \setminus \Per_A^n (X) = \sep_A^n (X)$, and $X \setminus \Per^n (X) = \sep^n (X)$. Furthermore $\sep_A^n (X) = \sep_A^{-n} (X)$ with similar equalities for $n$ and $-n$ ($n \in \mathbb{Z}$) holding for $\Per_A^n (X)$, $\sep^n (X)$ and $\Per^n (X)$ as well.
\begin{defn}\label{domun}
We say that a non-empty subset of $X$ is a \emph{domain of uniqueness for $A$} if every function in $A$ that vanishes on it, vanishes on the whole of $X$.
\end{defn}
For example, using results from elementary topology one easily shows that for a completely regular topological space $X$, a subset of $X$ is a domain of uniqueness for $C(X)$ if and only if it is dense in $X$.
\begin{thm}\label{commdesc}
The unique maximal abelian subalgebra of $A \rtimes_{\widetilde{\sigma}} \mathbb{Z}$ that contains $A$ is precisely the set of elements
\[A' =\{\sum_{n \in \mathbb{Z}} f_n \delta^n \,|\, \textup{for all} \, n \in \mathbb{Z}: {f_n}_{\upharpoonright \sep^n_A (X)} \equiv 0\}.\]
\end{thm}
\begin{proof}
Quoting a part of the proof of Proposition~\ref{commabel}, we have that $\sum_{n \in \mathbb{Z}} f_n \delta^n \in A'$ if and only if
\[\forall h \in A, \forall n \in \mathbb{Z} : f_n \cdot \widetilde{\sigma}^n (h) = f_n \cdot h.\]
Clearly this is equivalent to
\[\forall h \in A, \forall n \in \mathbb{Z} : {f_n}_{\upharpoonright \sep^n_A (X)} \equiv 0.\]
The result now follows from Proposition~\ref{commabel}.
\end{proof}
Note that for any non-zero integer $n$, the set $\{f_n \in A \, : \, {f_n}_{\upharpoonright \sep_A^n (X)} \equiv 0\}$ is a $\mathbb{Z}$-invariant ideal in $A$.
Note also that if $\sigma$ has finite order, then by Theorem~\ref{commdesc}, $A$ is not maximal abelian.
The following corollary follows directly from Theorem~\ref{commdesc}.
\begin{cor}\label{sepptdesc}
If $A$ separates the points of $X$, then $A'$ is precisely the set of elements
\[A' =\{\sum_{n \in \mathbb{Z}} f_n \delta^n \,|\, \textup{for all} \, n \in \mathbb{Z}: \supp(f_n) \subseteq \Per^n (X)\}.\]
\end{cor}
\begin{proof}
Immediate from Theorem~\ref{commdesc} and the remarks following Definition~\ref{sepper}.
\end{proof}
Applying Definition~\ref{domun} yields the following direct consequence of Theorem~\ref{commdesc}.
\begin{thm}\label{sepdom}
The subalgebra $A$ is maximal abelian in $A \rtimes_{\widetilde{\sigma}} \mathbb{Z}$ if and only if, for every $n \in \mathbb{Z} \setminus \{0\}$, $\sep_A^n (X)$ is a domain of uniqueness for $A$.
\end{thm}
In what follows we shall mainly focus on cases where $X$ is a topological space. Before we turn to such contexts, however, we use Theorem~\ref{commdesc} to give a description of the center of the crossed product, $Z (A \rtimes_{\widetilde{\sigma}} \mathbb{Z})$, for the general setup described in the beginning of this section.
\begin{thm}\label{center}
An element $g = \sum_{m \in  \mathbb{Z}} g_m \delta^m$ is in $Z(A \rtimes_{\widetilde{\sigma}} \mathbb{Z})$ if and only if both of the following conditions are satisfied:
\begin{enumerate}
\item for all $m \in \mathbb{Z}$, $g_m$ is $\mathbb{Z}$-invariant, and
\item for all $m \in \mathbb{Z}$, $ {g_m}_{\upharpoonright \sep^m_A (X)} \equiv 0$.
\end{enumerate}
\end{thm}
\begin{proof}
If $g$ is in $Z (A \rtimes_{\widetilde{\sigma}} \mathbb{Z})$ then certainly $g \in A'$, and hence condition (ii) follows from Theorem~\ref{commdesc}.
For condition (i), note that $g$ is in the center if and only if $g$ commutes with every element on the form $f_n \delta^n$.
Multiplying out, or looking at~(\ref{ekv}), we see that this means that
\begin{align*}&g = \sum_{m \in \mathbb{Z}} g_m \delta^m \,\,\textup{is in}\,\, Z(A \rtimes_{\widetilde{\sigma}} \mathbb{Z}) \Longleftrightarrow \\
&\forall n \in \mathbb{Z}, \forall m \in \mathbb{Z}, \forall f \in A: f \cdot \widetilde{\sigma}^n (g_m) = g_m \cdot \widetilde{\sigma}^m (f).
\end{align*}
We fix $m \in \mathbb{Z}$ and an $x \in \Per^m_A (X)$. Then for all $f \in A : f (x) = \widetilde{\sigma}^m (f)(x)$. If there is a function $f \in A$ that does not vanish in $x$, then for $g$ to be in the center we clearly must have that for all $n\in \mathbb{Z} : \,g_m (x) = \widetilde{\sigma}^n (g_m)(x)$. If all $f \in A$ vanish in $x$, then in particular both $g_m$ and $\widetilde{\sigma}^n (g_m)$ do. Thus for all points $x \in \Per^m_A (X)$ we have that $g_m$ is constant along the orbit of $x$ (i.e., for all $n \in \mathbb{Z} : g_m (x) = \widetilde{\sigma}^n (g_m) (x)$) for all $m \in \mathbb{Z}$, since $m$ was arbitrary in our above discussion. It remains to consider $x \in \sep^m_A (X)$. For such $x$, we have concluded that $g_m (x) =0$. If there exists $f \in A$ that does not vanish in $x$, we see that in order for the equality above to be satisfied we must have $\widetilde{\sigma}^n (g_m) (x) =0$ for all $n$, and if all $f \in A$ vanish in $x$, then in particular $\widetilde{\sigma}^n (g_m)$ does for all $n$ and the result follows.
\end{proof}
We now focus solely on topological contexts. The following theorem makes use of Corollary~\ref{sepptdesc}.
\begin{thm}\label{bairemax}
Let $X$ be a Baire space, and let $\sigma : X \rightarrow X$ be a homeomorphism inducing, as usual, an automorphism $\widetilde{\sigma}$ of $C(X)$.  Suppose $A$ is a subalgebra of $C(X)$ that is invariant under $\widetilde{\sigma}$ and its inverse, separates the points of $X$ and is such that for every non-empty open set $U \subseteq X$, there is a non-zero $f \in A$ that vanishes on the complement of $U$. Then $A$ is a maximal abelian subalgebra of  $A \rtimes_{\widetilde{\sigma}} \mathbb{Z}$ if and only if $\Per^{\infty} (X)$ is dense in $X$.
\end{thm}
\begin{proof}
Assume first that $\Per^{\infty} (X)$ is dense in $X$. This means in particular that any continuous function that vanishes on $\Per^{\infty} (X)$ vanishes on the whole of $X$. Thus Corollary~\ref{sepptdesc} tells us that $A$ is a maximal abelian subalgebra of $A \rtimes_{\widetilde{\sigma}} \mathbb{Z}$.
Now assume that $\Per^{\infty} (X)$ is \emph{not} dense in $X$. This means that $\bigcap_{n \in \mathbb{Z}_{>0}} (X \setminus \Per^n (X))$ is not dense. Note that the sets $X \setminus \Per^n(X)$, $n \in \mathbb{Z}_{>0}$ are all open. Since $X$ is a Baire space there exists an $n_0 \in \mathbb{Z}_{>0}$ such that $\Per^{n_0} (X)$ has non-empty interior, say $U \subseteq \Per^{n_0} (X)$. By the assumption on $A$, there is a nonzero function $f_{n_0} \in A$ that vanishes outside $U$. Hence Corollary~\ref{sepptdesc} shows that $A$ is \emph{not} maximal abelian in the crossed product.
\end{proof}
\begin{ex}\label{exloc}
Let $X$ be a locally compact Hausdorff space, and $\sigma : X \rightarrow X$ a homeomorphism. Then $X$ is a Baire space, and $C_c (X),\, C_0 (X),\, C_b (X)$ and $C(X)$ all satisfy the required conditions for $A$ in Theorem~\ref{bairemax}. For details, see for example \cite{munkres}. Hence these function algebras are maximal abelian in their respective crossed products with $\mathbb{Z}$ under $\sigma$ if and only if $\Per^{\infty} (X)$ is dense in $X$.
\end{ex}
\begin{ex}\label{irrot}
Let $X = \mathbb{T}$ be the unit circle in the complex plane, and let $\sigma$ be counterclockwise rotation by an angle which is an irrational multiple of $2 \pi$. Then every point is non-periodic and thus, by Theorem~\ref{bairemax}, $C(\mathbb{T})$ is maximal abelian in the associated crossed product.
\end{ex}
\begin{ex}\label{ratrot}
Let $X=\mathbb{T}$ and $\sigma$ counterclockwise rotation by an angle which is a rational multiple of $2 \pi$, say $2 \pi p/q$ (where $p, q$ are relatively prime positive integers). Then every point on the circle has period precisely $q$, and the non-periodic points are certainly not dense. Using Corollary~\ref{sepptdesc} we see that
\[C(\mathbb{T})' =\{\sum_{n \in I} f_{nq} \delta^{nq} \,|\, f_{nq} \in C(\mathbb{T})\}.\]
\end{ex}
We will use the following theorem to display an example different in nature from the ones already considered.
\begin{thm}\label{domunmax}
Let $X$ be a topological space, $\sigma : X \rightarrow X$ a homeomorphism, and $A$ a non-zero subalgebra of $C(X)$, invariant both under the usual induced automorphism $\widetilde{\sigma}: C(X) \rightarrow C(X)$ and under its inverse. Assume that $A$ separates the points of $X$ and is such that every non-empty open set $U \subseteq X$ is a domain of uniqueness for $A$. Then $A$ is maximal abelian in $A \rtimes_{\widetilde{\sigma}} \mathbb{Z}$ if and only if $\sigma$ is not of finite order (that is, $\sigma^n \neq \id_X$ for any non-zero integer $n$).
\end{thm}
\begin{proof}
By Corollary~\ref{sepptdesc}, $A$ being maximal abelian implies that $\sigma$ is not of finite order. Indeed, if $\sigma^p = \id_X$, where $p$ is the smallest such positive integer, then $X = \Per^p (X)$ and hence $f \delta^p \in A'$ for any $f \in A$.
For the converse, assume that $\sigma$ does not have finite order. The sets $\sep^n (X)$ are non-empty open subsets of $X$ for all $n \neq 0$ and thus domains of uniqueness for $A$ by assumption of the theorem. The implication now follows directly from Corollary~\ref{sepptdesc}.
\end{proof}
\begin{cor}\label{comfd}
Let $M$ be a connected complex manifold and suppose the function $\sigma : M \rightarrow M$ is biholomorphic. If $A \subseteq H(M)$ is a subalgebra of the algebra of holomorphic functions that separates the points of $M$ and which is invariant under the induced automorphism $\widetilde{\sigma}$ of $C(M)$ and its inverse, then $A \subseteq A \rtimes_{\widetilde{\sigma}} \mathbb{Z}$ is maximal abelian if and only if $\sigma$ is not of finite order.
\end{cor}
\begin{proof}
On connected complex manifolds, open sets are domains of uniqueness for $H(M)$. See for example \cite{complex}.
\end{proof}
\begin{rmk}\label{disrmk}
It is important to point out that the required conditions in Theorem~\ref{bairemax} and Theorem~\ref{domunmax} can only be simultaneously satisfied in case $X$ consists of a single point and $A = \mathbb{C}$.
To see this we assume that both conditions are satisfied. Note first of all that this implies that every non-empty open subset of $X$ is dense. Assume to the contrary that there is a non-empty open subset $U \subseteq X$ that is not dense in $X$. We may then choose a non-zero $f \in A$ that vanishes on $X \setminus U$. Certainly, $f$ must then vanish on $V = X \setminus \overline{U}$. As $U$ is not dense, however, this implies that $f$ is identically zero since the non-empty open set $V$ is a domain of uniqueness by assumption. Hence we have a contradiction and conclude that every non-empty open subset of $X$ is dense.
Secondly, we note that since $\mathbb{C}$ is Hausdorff and $A \subseteq C(X)$ separates the points of $X$, $X$ must be Hausdorff.
Assume now that there are two distinct points $p, q \in X$. As $X$ is Hausdorff we can separate them by two disjoint open sets. Since every non-empty open subset is dense, however, this is not possible and hence $X$ consists of one point. If $X = \{p\}$, the only possibilities for $A$ are $\{0\}$ and $\mathbb{C}$. As the conditions in Theorem~\ref{bairemax} imply existence of non-zero functions in $A$ we conclude that $A = \mathbb{C}$.
Conversely, if $X = \{p\}$ and $A = \mathbb{C}$, the conditions in both Theorem~\ref{bairemax} and Theorem~\ref{domunmax} are satisfied.
\end{rmk}
\section{Automorphisms of commutative completely regular semi-simple Banach algebras}
\subsection{Motivation} \label{sec:motivation}
In the setup in Theorem~\ref{bairemax} we concluded that we have an appealing equivalence between density of the non-periodic points in $X$ and $A$ being maximal abelian in the associated crossed product, under a certain condition on $A$. Example~\ref{irrot} and Example~\ref{ratrot}, respectively, show an instance where both these equivalent statements are true, and where they are false.
Generally speaking, however, the density of the non-periodic points and $A$ being maximal abelian are unrelated properties; all four logical possibilities can occur. We show this by giving two additional examples.
\begin{ex}\label{smallalg}
As in Example~\ref{irrot}, let $X = \mathbb{T}$ be the unit circle and $\sigma$ counterclockwise rotation by an angle that is an irrational multiple of $2 \pi$. If we use $A = \mathbb{C}$ instead of $C(\mathbb{T})$, the whole crossed product is commutative and thus $A$ is clearly \emph{not} maximal abelian in it. The non-periodic points, however, are of course still dense. Here we simply chose a subalgebra of $C(\mathbb{T})$ so small that the homeomorphism $\sigma$ was no longer visible in the crossed product.
\end{ex}
\begin{ex}\label{addpt}
Let $X = \mathbb{T} \cup {\{0\}}$ with the usual subspace topology from $\mathbb{C}$, and let $\sigma'$ be such that it fixes the origin and rotates points on the circle counterclockwise with an angle that is an irrational multiple of $2 \pi$. As function algebra $A$, we take $C(\mathbb{T})$ and extend every function in it to $X$ so that it vanishes in the origin. This is obviously an algebra of functions being continuous on $X$, which is invariant under $\sigma'$ and its inverse. Since $A$ separates points in $X$, Corollary~\ref{sepptdesc} assures us that $A$ is maximal abelian in the crossed product, even though the non-periodic points are not dense in $X$.
\end{ex}
In the following example, the equivalence in Theorem~\ref{bairemax} fails in the same fashion as in Example~\ref{addpt}. It is included, however, since it will be illuminating to refer to it in what follows.
\begin{ex}\label{holom}
As in Example~\ref{addpt}, let $X = \mathbb{T} \cup \{0\}$ and
$\sigma'$ the map defined as counterclockwise rotation by an angle that is an irrational multiple of $2 \pi$ on $\mathbb{T}$ and $\sigma' (0) = 0$. Let $A$ be the
restriction to $X$ of all continuous functions on the closed unit disc
that are holomorphic on the open unit disc. By the maximum modulus theorem,
none of these functions are non-zero solely in the origin. Thus, by
Corollary~\ref{sepptdesc}, we again obtain a case where the non-periodic points
are not dense, but $A$ is a maximal abelian subalgebra in the crossed product.
\end{ex}
In summary, we have now displayed three examples where we do not have an equivalence between algebra and topological dynamics as in Theorem~\ref{bairemax}. In the following subsection we prove a general result - in the context of automorphisms of Banach algebras - that in particular shows that for a certain class of pairs of discrete dynamical systems and $\mathbb{Z}$-invariant function algebras on it, $((X, \sigma),  A)$, yielding the associated crossed product $A \rtimes_{\widetilde{\sigma}} \mathbb{Z}$ as usual, one can always find another such pair $((Y, \phi),B)$ with associated crossed product $B \rtimes_{\widetilde{\phi}} \mathbb{Z}$ \emph{canonically isomorphic} to $A \rtimes_{\widetilde{\sigma}} \mathbb{Z}$, where the equivalence \emph{does} hold: the non-periodic points of $Y$ are dense if and only if $B$ is maximal abelian in $B \rtimes_{\widetilde{\phi}} \mathbb{Z}$ (which it is, by the canonical isomorphism, if and only if $A$ is maximal abelian in $A \rtimes_{\widetilde{\sigma}} \mathbb{Z}$). In this way, the equivalence of an algebraic property with a topological dynamical property is restored. Examples~\ref{smallalg} through~\ref{holom} all fall into this mentioned class of pairs, as we will see in Example~\ref{dualc} -~\ref{dualholom}.
\subsection{A system on the character space}
\label{sec:dualsystem}
We will now focus solely on Banach algebras, and start by recalling a number of basic results concerning them. We refer to \cite{larsen} for details. All Banach algebras under consideration will be \emph{complex} and \emph{commutative}.
\begin{defn}\label{maxid}
Let $A$ be a complex commutative Banach algebra. The set of all non-zero multiplicative linear functionals on $A$ is denoted by $\Delta(A)$ and called the \emph{character space} of $A$.
\end{defn}
\begin{defn}\label{gelfand}
Given any $a \in A$, we define a function $\widehat{a} : \Delta(A) \rightarrow \mathbb{C}$ by $\widehat{a} (\mu) = \mu(a) \,\, (\mu \in \Delta(A))$. The function $\widehat{a}$ is called the \emph{Gelfand transform of $a$}.
Let $\widehat{A} = \{\widehat{a} \,|\, a \in A\}$. The character space $\Delta(A)$ is endowed with the topology generated by $\widehat{A}$, which is called the \emph{Gelfand topology} on $\Delta(A)$. The Gelfand topology is locally compact and Hausdorff. A commutative Banach algebra $A$ for which the Gelfand transform, i.e., the map sending $a$ to $\widehat{a}$, is injective, is called \emph{semi-simple}.
\end{defn}
Let $A$ be a commutative semi-simple complex Banach algebra and $\widetilde{\sigma}: A \rightarrow A$ an algebra automorphism. Then $\widetilde{\sigma}$ induces a bijection $\sigma : \Delta(A) \rightarrow \Delta(A)$ defined by $\sigma (\mu) = \mu \circ \widetilde{\sigma}^{-1}, \,\, (\mu \in \Delta(A))$, which is automatically a homeomorphism when $\Delta(A)$ has the Gelfand topology.
Note that by semi-simplicity of $A$, the map \[\phi : \Aut (A) \rightarrow \{\sigma \in \Homeo(\Delta(A)) \,|\, \widehat{a} \circ \sigma,\, \widehat{a} \circ \sigma^{-1} \in \widehat{A} \,\,\textup{for all}\,\,  a \in A\}\] defined by \[\phi(\widetilde{\sigma}) (\mu) = \mu \circ \widetilde{\sigma}^{-1}\] is an isomorphism of groups.
In turn, $\sigma$ induces an automorphism $\widehat{\sigma}$ on $\widehat{A}$ as in Section~\ref{auto}, namely $\widehat{\sigma} (\widehat{a}) = \widehat{a} \circ \sigma^{-1} = \widehat{\widetilde{\sigma} (a)}$.

The following result shows that in the context of a semi-simple Banach algebra one may pass to an isomorphic crossed product, but now with an algebra of continuous functions on a topological space. It is here that topological dynamics can be brought into play again. The proof consists of a trivial direct verification.
\begin{thm}\label{isom}
Let $A$ be a commutative semi-simple Banach algebra and $\widetilde{\sigma}$ an automorphism, inducing an automorphism $\widehat{\sigma} : \widehat{A} \rightarrow \widehat{A}$ as above. Then the map $\Phi : A \rtimes_{\widetilde{\sigma}} \mathbb{Z} \rightarrow \widehat{A} \rtimes_{\widehat{\sigma}} \mathbb{Z}$ defined by $\sum_{n \in \mathbb{Z}} a_n \delta^n \mapsto \sum_{n \in \mathbb{Z}} \widehat{a_n} \delta^n$ is an isomorphism of algebras mapping $A$ onto $\widehat{A}$.
\end{thm}
\begin{defn}\label{compreg}
A commutative Banach algebra $A$ is said to be \emph{completely regular} if for every subset $F \subseteq \Delta(A)$ that is closed in the Gelfand-topology and for every $\phi_0 \in \Delta(A) \setminus F$ there exists an $a \in A$ such that $\widehat{a} (\phi) = 0$ for all $\phi\in F$ and $\widehat{a} (\phi_0) \neq 0$. In Banach algebra theory it is proved that $A$ is completely regular if and only if the hull-kernel topology on $\Delta(A)$ coincides with the Gelfand topology, see for example \cite{bondun}.
\end{defn}
\begin{thm}\label{compregmax}
Let $A$ be a commutative completely regular semi-simple Banach algebra, $\widetilde{\sigma}: A \rightarrow A$ an algebra automorphism and $\sigma$ the homeomorphism on $\Delta(A)$ in the Gelfand topology induced by $\widetilde{\sigma}$ as described above. Then the non-periodic points of $(\Delta(A), \sigma)$ are dense if and only if $\widehat{A}$ is a maximal abelian subalgebra of $\widehat{A} \rtimes_{\widehat{\sigma}} \mathbb{Z}$. In particular, $A$ is maximal abelian in $A \rtimes_{\widetilde{\sigma}} \mathbb{Z}$ if and only if the non-periodic points of $(\Delta(A), \sigma)$ are dense.
\end{thm}
\begin{proof}
As mentioned in Definition~\ref{gelfand}, $\Delta(A)$ is locally compact Hausdorff in the Gelfand topology, and clearly $\widehat{A}$ is by definition a separating function algebra on it. Since we assumed $A$ to be completely regular, we see that all the conditions assumed in Theorem~\ref{bairemax} are satisfied, and thus this theorem yields the equivalence. Furthermore, by Theorem~\ref{isom}, $A$ is maximal abelian in $A \rtimes_{\widetilde{\sigma}} \mathbb{Z}$ if and only if $\widehat{A}$ is maximal abelian in $\widehat{A} \rtimes_{\widehat{\sigma}} \mathbb{Z}$.
\end{proof}
We shall now revisit Examples~\ref{smallalg} through~\ref{holom} and use Theorem~\ref{compregmax} to conclude algebraic properties from topological dynamics after all.
\begin{ex}\label{dualc}
Consider again Example~\ref{smallalg}. Obviously $\Delta(\mathbb{C}) = \{\id_{\mathbb{C}}\}$. Thus $\mathbb{C}$ is a
commutative semi-simple completely regular Banach algebra. Trivially, a set
with only one element has no non-periodic point, so that $A = \mathbb{C}$ is
not maximal abelian by Theorem~\ref{compregmax}.
\end{ex}
\begin{ex}\label{dualaddpt}
Consider again Example~\ref{addpt}. Clearly the function algebra $A$ on $X$ is isometrically isomorphic to $C(\mathbb{T})$ (when both algebras are endowed with
the sup-norm) and thus a commutative completely regular semi-simple Banach
algebra. It is furthermore a well known result from the theory of Banach
algebras that $\Delta(C(X)) = \{\mu_x \,|\, x \in X\}$ for any compact
Hausdorff space $X$, where $\mu_x$ denotes the point evaluation in $x$, and
that $\Delta(C(X))$ equipped with the Gelfand topology is homeomorphic to
$X$ (see for example \cite{larsen}). Thus clearly here $\Delta(A) =
\Delta(C(\mathbb{T})) = \mathbb{T}$. Of course the induced map $\sigma$
defined on $\Delta(A)$ by $\sigma (\mu_x) = \mu_{\sigma(x)}$ can 
simply be identified with rotation $\sigma'$ of the circle by an angle that is an irrational multiple of $2 \pi$. So on
$\Delta(A)$ the non-periodic points are dense, and hence $A$ is maximal abelian
by Theorem~\ref{compregmax}. Here passing to the system on the character space
corresponded to deleting the origin from $X$, thus recovering $\mathbb{T}$, and
restricting $A$ to $\mathbb{T}$, hence recovering $C(\mathbb{T})$ so that in
the end we recovered the setup in Example~\ref{irrot}.
\end{ex}
\begin{ex}\label{dualholom}
Consider again Example~\ref{holom}. Using the maximum modulus theorem one sees
that $A$ is isometrically isomorphic to $\mathcal{A}(\overline{\mathbb{D}})$, the algebra of all functions
that are continuous on the unit circle and holomorphic on the open unit disc (denoted by $\mathbb{D}$), as $A$
is the restrictions of such functions to $\mathbb{T} \cup \{0\}$. Hence $A$ is
a commutative completely regular semi-simple  Banach algebra. Furthermore it is
a standard result from Banach algebra theory that $\Delta(\mathcal{A}(\overline{\mathbb{D}}))$ (endowed with
the Gelfand topology) is canonically homeomorphic to $\overline{\mathbb{D}}$; the elements in
$\Delta(\mathcal{A}(\overline{\mathbb{D}}))$ are precisely the point evaluations in $\overline{\mathbb{D}}$ (see \cite{larsen}).
So we conclude that $\Delta(A)$ is also equal to $\overline{\mathbb{D}}$. The induced
homeomorphism $\sigma$ on $\Delta(A) = \overline{\mathbb{D}}$ is rotation by the same angle as for 
$\sigma'$. Here the non-periodic points are obviously dense, so that $A$ is
maximal abelian by Theorem~\ref{compregmax}.
\end{ex}
Note the difference in nature between Examples~\ref{dualaddpt} and ~\ref{dualholom}. In the former, passing to the system on the character space corresponds to deleting a point from the original system and restricting the function algebra and homeomorphism, while in the latter it corresponds to adding (a lot of) points and extending.

We conclude this subsection by giving yet another example of an application of Theorem~\ref{compregmax}, recovering one of the results we obtained in Example~\ref{exloc} by using Theorem~\ref{bairemax}.
\begin{ex}\label{locc0}
Let $X$ be a locally compact Hausdorff space, and $\sigma : X \rightarrow X$ a homeomorphism. Let $A = C_0 (X)$; then $\sigma$ induces an automorphism $\widetilde{\sigma}$ on $A$ as in Section~\ref{auto}. Here $\Delta(A)$  is canonically homeomorphic $X$ and $A$ is canonically isomorphic (with respect to the homeomorphism between $\Delta(A)$ and $X$) to $\widehat{A}$ (see for example \cite{larsen}). Hence by Theorem~\ref{compregmax} $A$ is maximal abelian in $A \rtimes_{\widetilde{\sigma}} \mathbb{Z}$ if and only if the non-periodic points of $X$ are dense, as already mentioned in Example~\ref{exloc}.
\end{ex}
\subsection{Integrable functions on locally compact abelian groups}
\label{sec:loccompgr}
In this subsection we consider the crossed product $L_1 (G) \rtimes_{\Psi} \mathbb{Z}$, where $G$ is a locally compact abelian group and $\Psi : L_1 (G) \rightarrow L_1 (G)$ an automorphism. We will show that under an additional condition on $G$, a stronger result than Theorem~\ref{compregmax} is true (cf. Theorem~\ref{conndual}).

We start by recalling a number of standard results from the theory of Fourier
analysis on groups, and refer to \cite{larsen} and \cite{rudin} for details.
Let $G$ be a locally compact abelian group. Recall that $L_1 (G)$ consists of
equivalence classes of complex valued Borel measurable functions of $G$ that
are integrable with respect to a Haar measure on $G$, and that $L_1 (G)$
equipped with convolution product is a commutative completely regular
semi-simple Banach algebra. A group homomorphism $\gamma : G \rightarrow \mathbb{T}$ from a locally compact
abelian group $G$ to the unit circle is called a \emph{character} of $G$.
The set of all \emph{continuous} characters of $G$ forms a group $\Gamma$, the \emph{dual group} of $G$, if the group operation is defined by
 \[(\gamma_1 + \gamma_2) (x) = \gamma_1 (x) \gamma_2 (x) \quad (x \in G; \,\gamma_1, \gamma_2 \in \Gamma).\]
If $\gamma \in \Gamma$ and if
\[\widehat{f}(\gamma) = \int_{G} f(x) \gamma (-x) dx \quad (f \in L_1 (G)),\]
then the map $f \mapsto \widehat{f} (\gamma)$ is a non-zero complex homomorphism of $L_1 (G)$. Conversely, every non-zero complex homomorphism of $L_1 (G)$ is obtained in this way, and distinct characters induce distinct homomorpisms. Thus we may identify $\Gamma$ with $\Delta(L_1 (G))$.
The function $\widehat{f} : \Gamma \rightarrow \mathbb{C}$ defined as above is called the \emph{Fourier transform} of $f, \, f \in L_1 (G)$, and is hence precisely the Gelfand transform of $f$. We denote the set of all such $\widehat{f}$ by $A(\Gamma)$. Furthermore, $\Gamma$ is a locally compact abelian group in the Gelfand topology.
\begin{defn}\label{ring}
Given a set $X$, a \emph{ring of subsets} of $X$ is a collection of subsets of $X$ which is closed under the formation of finite unions, finite intersections, and complements (in $X$). Note that any intersection of rings is again a ring.
The \emph{coset-ring} of $\Gamma$ is defined to be the smallest ring of subsets of $\Gamma$ which contains all \emph{open} cosets, i.e., all subsets of $\Gamma$ of the form $a + U$, where $a \in \Gamma$ and $U$ is an open subgroup of $\Gamma$.
\end{defn}
We are now ready to define a particular type of map on the coset ring of $\Gamma$ (cf. \cite{rudin}).
\begin{defn}\label{pwaffine}
Let $E$ be a coset in $\Gamma$. A continuous map $\sigma : E \rightarrow \Gamma$ which satisfies the identity
\[\sigma (\gamma + \gamma' - \gamma'') = \sigma(\gamma) + \sigma(\gamma') - \sigma(\gamma'') \quad (\gamma, \gamma', \gamma'' \in E)\]
is called \emph{affine}.
Suppose that
\begin{enumerate}
\item $S_1, \ldots, S_n$ are pairwise disjoint sets belonging to the coset-ring of $\Gamma$;
\item each $S_i$ is contained in an open coset $K_i$ in $\Gamma$;
\item for each $i$, $\sigma_i$ is an affine map of $K_i$ into $\Gamma$;
\item $\sigma$ is the map of $Y = S_1 \cup \ldots \cup S_n$ into $\Gamma$ which coincides on $S_i$ with $\sigma_i$.
\end{enumerate}
Then $\sigma$ is said to be a \emph{piecewise affine} map from $Y$ to $\Gamma$.
\end{defn}
The following theorem is a key result for what follows. It states that every automorphism of $L_1 (G)$ is induced by a piecewise affine homeomorphism, and that a piecewise affine homeomorphism induces an injective homomorphism from $L_1 (G)$ to itself.
\begin{thm}\label{autohomeo}
Let $\widetilde{\sigma} : L_1 (G) \rightarrow L_1 (G)$ be an automorphism. Then for every $f \in L_1 (G)$ we have that $\widehat{\widetilde{\sigma} (f)} = \widehat{f} \circ \sigma$, where $\sigma : \Gamma \rightarrow \Gamma$ is a fixed piecewise affine homeomorphism. Also, if $\sigma : \Gamma \rightarrow \Gamma$ is a piecewise affine homeomorphism, then $\widehat{f} \circ \sigma \in A(\Gamma)$ for every $\widehat{f} \in A(\Gamma)$.
\end{thm}
\begin{proof}
Follows from the more general Theorems 4.1.3 and 4.6.2 in \cite{rudin}.
\end{proof}
Now let $\widetilde{\sigma} : L_1 (G) \rightarrow L_1(G)$ be an automorphism
and consider the crossed product $L_1 (G) \rtimes_{\widetilde{\sigma}}
\mathbb{Z}$. Letting $\widetilde{\sigma}$ induce a homeomorphism as described
for arbitrary commutative completely regular semi-simple Banach algebras in the
paragraph following Definition~\ref{gelfand}, we obtain
$\sigma^{-1} : \Gamma \rightarrow \Gamma$, where $\sigma$ is the piecewise
affine homeomorphism inducing $\widetilde{\sigma}$ in accordance with
Theorem~\ref{autohomeo}.

\begin{thm}\label{conndual}
Let $G$ be a locally compact abelian group with connected dual group and let $\widetilde{\sigma} : L_1 (G) \rightarrow L_1 (G)$ be an automorphism. Then $L_1 (G)$ is maximal abelian in $L_1 (G) \rtimes_{\widetilde{\sigma}} \mathbb{Z}$ if and only if $\widetilde{\sigma}$ is not of finite order.\end{thm}
\begin{proof}
Denote by $\Gamma$ the dual group of $G$.
By Theorem~\ref{compregmax}, $\sigma$, the homeomorphism induced by $\widetilde{\sigma}$ in accordance with the discussion following Definition~\ref{gelfand}, is not of finite order if $L_1 (G)$ is maximal abelian. Assume now that $L_1 (G)$ is not maximal abelian. By Theorem~\ref{compregmax}, this implies that $\Per^{\infty} (\Gamma)$ is not dense in $\Gamma$. As argued in the proof of Theorem~\ref{bairemax}, there must then exist $n_0 \in \mathbb{N}$ such that $\Per^{n_0} (\Gamma)$ has non-empty interior. Namely, since in this case $\bigcap_{n \in \mathbb{Z}_{>0}} (\Gamma \setminus \Per^n (\Gamma))$ is not dense and that the sets $\Gamma \setminus \Per^n(\Gamma)$, $n \in \mathbb{Z}_{>0}$ are all open, the fact that $\Gamma$ is a Baire space (being locally compact and Hausdorff) implies existence of an $n_0 \in \mathbb{Z}_{>0}$ such that $\Per^{n_0} (\Gamma)$ has non-empty interior. Note that $\Gamma$ being connected implies that $\sigma^{-1}$, the piecewise affine homeomorphism of $\Gamma$ inducing $\widetilde{\sigma}$ in accordance with Theorem~\ref{autohomeo}, must be affine by connectedness of $\Gamma$ (the coset-ring is trivially $\{\emptyset, \Gamma\}$) and hence so is $\sigma$. It is readily verified that the map $\sigma^{n_0} - I$ is then also affine. Now clearly $\Per^{n_0} (\Gamma) = (\sigma^{n_0} - I )^{-1} (\{0\})$.
The affine nature of $\sigma^{n_0} - I$ assures us that $\Per^{n_0} (\Gamma)$ is a coset. In a topological group, however, continuity of the group operations implies that cosets with non-empty interior are open, hence also closed. We conclude that $\Per^{n_0} (\Gamma)$ is a non-empty closed and open set. Connectedness of $\Gamma$ now implies that every point in $\Gamma$ is $n_0$-periodic under $\sigma$. Hence, by the discussion following Definition~\ref{gelfand}, ${\widetilde{\sigma}}^{n_0}$ is the identity map on $L_1 (G)$.
\end{proof}
The following example shows that if the dual of $G$ is not connected, the equivalence in Theorem~\ref{conndual} need not hold.
\begin{ex}\label{disco}
Let $G = \mathbb{T}$ be the circle group. Here $\Gamma = \mathbb{Z}$ (see \cite{larsen} for details), which is not connected.
Define $\sigma : \mathbb{Z} \rightarrow \mathbb{Z}$ by $\sigma (n) = n \,\, (n \in 2 \mathbb{Z})$ and $\sigma (m) = m+2\,\, (m \in 1 + 2 \mathbb{Z})$. Obviously $\sigma$ and $\sigma^{-1}$ are then piecewise affine homeomorphisms that are not of finite order. By Theorem~\ref{autohomeo}, $\sigma$ induces an automorphism $\widetilde{\sigma^{-1}} : L_1 (\mathbb{T}) \rightarrow L_1 (\mathbb{T})$, which in turn induces the homeomorphism $\sigma^{-1} : \mathbb{Z} \rightarrow \mathbb{Z}$. Now $A(\mathbb{Z})$ is not maximal abelian in $A(\mathbb{Z}) \rtimes_{\widehat{\sigma^{-1}}} \mathbb{Z}$ since by Corollary~\ref{sepptdesc} we have
\[A(\mathbb{Z})' =\{\sum_{n \in \mathbb{Z}} \widehat{f_n} \delta^n \,|\, \textup{for all} \, n \in \mathbb{Z} \setminus \{0\}: \supp(\widehat{f_n}) \subseteq 2 \mathbb{Z}\},\]
and hence by Theorem~\ref{isom}
\[L_1 (\mathbb{T})' =\{\sum_{n \in \mathbb{Z}} f_n \delta^n \,|\, \textup{for all} \, n \in \mathbb{Z} \setminus \{0\}: \supp(\widehat{f_n}) \subseteq 2 \mathbb{Z}\}.\]
Note that
\[\{\sum_{n \in \mathbb{Z}} f_n \delta^n \,|\, \textup{for all} \, n \in \mathbb{Z} \setminus \{0\}: f_n \in \mathbb{C} [z^2, z^{-2}]\} \subseteq L_1 (\mathbb{T})',\]
and thus we conclude that $L_1 (\mathbb{T})$ is not maximal abelian.
\end{ex}
\subsection{A theorem on generators for the commutant}
\label{sec:generatorscom}
A natural question to ask is whether or not the commutant $A' \subseteq A \rtimes_{\widetilde{\sigma}} \mathbb{Z}$ is finitely generated as an algebra over $\mathbb{C}$ or not. Here we give an answer in the case when $A$ is a semi-simple commutative Banach algebra.
\begin{thm}\label{generation}
Let $A$ be a semi-simple commutative Banach algebra and let $\widetilde{\sigma}
: A \rightarrow A$ be an automorphism. Then $A'$ is finitely generated as an algebra
over $\mathbb{C}$ if and only if $A$ has finite dimension as a vector space
over $\mathbb{C}$.
\end{thm}
\begin{proof}
Let the induced homeomorphism $\sigma : \Delta(A)\rightarrow \Delta(A)$ be as
usual. Assume first that $A$ has infinite dimension. By basic theory of Banach
spaces, $A$ must then have uncountable dimension. If $A'$ were generated by
finitely many elements, then $A'$, and in particular $A$, would have countable
dimension, which is a contradiction. Hence $A'$ is not finitely generated as an
algebra over $\mathbb{C}$. For the converse, we need two results from Banach
algebra theory. Suppose that $A$ has finite dimension. By Proposition 26.7 in
\cite{bondun}, $\Delta(A)$ is then a finite set, and thus every point in
$(\Delta(A), \sigma)$ has a common finite period, $n_0$ say. Furthermore, by
Corollary 21.6 in \cite{bondun} $A$ must then also be unital. We pass now to
the crossed product $\widehat{A} \rtimes_{\widehat{\sigma}} \mathbb{Z}$, which
is isomorphic to $A \rtimes_{\widetilde{\sigma}} \mathbb{Z}$ by
Theorem~\ref{isom}. Clearly $\widehat{A}$ is unital and has finite linear
dimension since $A$ does. By Corollary~\ref{sepptdesc}, for a general element
$\sum_{n \in \mathbb{Z}} \widehat{a_n} \delta^n \in (\widehat{A})'$ the set of
possible coefficients of $\delta^n$ is a vector subspace (and even an ideal) of
$\widehat{A}$, $K_n$ say, and hence of finite dimension. Since all elements of
$(\Delta(A), \sigma)$ have period $n_0$, Corollary~\ref{sepptdesc} also tells
us that $K_{r + l \cdot n_0} = K_{r}$ for all $r, l \in \mathbb{Z}$. Now note
that since $\widehat{A}$ is unital, $\delta^{n_0}, \delta^{-n_0} \in
(\widehat{A})'$. Thus, denoting a basis for a $K_l$ by $\{e_{(l,1)}, \ldots,
e_{(l,l_r)}\}$ (where $l_r \leq s$), the above reasoning assures us that
$\bigcup_{l = 1}^{n_0} \bigcup_{j=1}^{l_r} \{e_{(l,j)} \delta^{l}\}$ generates
$(\widehat{A})'$ as an algebra over $\mathbb{C}$. By Theorem~\ref{isom} this
implies that also $A' \subseteq A \rtimes_{\widetilde{\sigma}} \mathbb{Z}$ is
finitely generated as an algebra over $\mathbb{C}$.
\end{proof}
\section*{Acknowledgments}
This work was supported by a visitor's grant of the \textit{Netherlands Organisation for Scientific Research (NWO)}, \textit{The Swedish Foundation for International Cooperation in Research and Higher Education (STINT)}, \textit{Crafoord Foundation and The Royal Physiographic Society in Lund}.

\end{document}